\RequirePackage{fix-cm}
\documentclass[smallextended]{svjour3}       
\smartqed  
\usepackage{amsmath}
\usepackage[utf8]{inputenc}
\usepackage{hyperref}

\hypersetup{
	colorlinks=true,
	linkcolor=blue,
	filecolor=blue,
	urlcolor=blue,
	citecolor=blue,
}
\usepackage{booktabs}
\usepackage[pdftex]{graphics}
\graphicspath{{figures}}
\usepackage{amsfonts}

\usepackage{epstopdf}
\usepackage[section]{placeins}
\usepackage{adjustbox}
\usepackage[a4paper, total={6in, 8in}]{geometry}
\usepackage{cite}
\usepackage{breqn}
\usepackage{framed}
\usepackage{caption}
\usepackage{epstopdf}
\usepackage{subcaption}
\usepackage{lipsum}
\usepackage{booktabs}
\usepackage{longtable}
\begin{document}

\title{ A Two-stage Network Data Envelopment Analysis: An Education Sector Application}

\author{Awadh Pratap Singh         \and
        Shiv Prasad Yadav
}

\institute{ Awadh Pratap Singh\at
              Department of Mathematics, Indian Institute of Technology Roorkee, Roorkee-247667, India. \\
              \email{awadhma2015@gmail.com}          
           \and
           Shiv Prasad Yadav \at
              Department of Mathematics, Indian Institute of Technology Roorkee, Roorkee-247667, India.\\
              \email{spyorfma@gmail.com} 
}

\date{Received: date / Accepted: date}

\maketitle

\begin{abstract}
 In general, data envelopment analysis (DEA) works like a black box that does not provide any adequate detail to identify the specific reason for inefficiency in decision-making units (DMUs). The motivation of this study is to analyze the cause of the inefficiency of DMUs in a decision process with the help of a two-stage relational network DEA  model. In the current study, a two-stage relational network DEA model is applied to measure the performance of DMUs for the whole process and each stage independently. In general, past studies used conventional DEA models in education sectors to analyze the performances of educational institutions. In the current study, by considering quantitative attributes to measure the performance of Indian institutes of management (IIMs) by using network DEA, we develop a procedure that captures both quality and quantity. 
 \\
 
\end{abstract}

\section{Introduction} \label{itm:intro}
Data envelopment analysis (DEA) is a data-supported nonparametric approach that is used to measure the relative performance of decision-making units (DMUs) \cite{cook2009data}.  Charnes, Cooper, and Rhodes \cite{charnes1978measuring} developed a model called the CCR DEA Model to measure the performance of DMUs.  The efficiency of DMU is defined as the ratio of output to input (efficiency = output/input). The relative efficiency of a DMU is defined as, the ratio of its efficiency to the largest efficiency under consideration. It lies in the interval (0,1]. A DMU with an efficiency score of 1 is called efficient DMU otherwise inefficient DMU.

Education and research play a vital role in the development of any nation. In the current study, we will try to measure the performance efficiencies of IIMs with the help of the network DEA model, particularly the two-stage network DEA. The main reason for using the two-stage network DEA model over the conventional CCR DEA model is to understand the causes behind the inefficiencies of DMUs. In general, DEA works like a black box that does not provide any adequate detail to identify the specific reason for inefficiency in DMUs. The motivation of this study is to analyze the cause of the inefficiency of DMUs in a decision process with the help of a two-stage relational network DEA model. There are many types of research available in literature where the network DEA model is developed and used for the performance evaluation of DMUs. Kao \cite{kao2014network} reported a detailed review on network DEA which consists of the DMUs with series, parallel and mixed structure.  Application of network DEA approach can be found in several fields like education sector, banking sector, healthcare sector, hotel industries, etc. Tan and Despotis \cite{tan2021investigation} investigated the efficiency of the hotel industry by using a network DEA model. Kao and Hung \cite{kao2008efficiency}, Arya and Singh \cite{arya2021development}, and Ma and Chen \cite{ma2018evaluating} proposed a two-stage parallel series network and provided the application of proposed models in the healthcare sector. Puri et al. \cite{puri2017new}, Lofti and Ghasemi \cite{lotfi2013multi} also used multi-component DEA in banking sector. 

In past studies, many researchers \cite{singh2021performance,arcelus1997efficiency,sinuany1994academic,bessent1982application,tomkins1988experiment,tyagi2009relative} used conventional DEA in education sector. As Lee and Worthington \cite{lee2016network} explained in their work that higher educational institutions are the more complex structure that uses multiple inputs to produce multiple outputs. The research process is also complicated. As, in some stage, some output of research process become the inputs for the next stages to produce final outputs. Many researchers applied network DEA models for performance evaluation in the education sector. Lee and Worthington \cite{lee2016network} applied a  network DEA model for evaluating the research performance of Australian universities. In this study, they developed a university research production model. The 1st stage uses input (full-time equivalent (FTE) academic staff and capital stock) to produce output as the number of publications. In the 2nd stage, the number of publications is used to produce the final output (grant secured). Our study is also motivated by Lee and Worthington's study. Tavares et al. \cite{tavares2021proposed} proposed a multistage network DEA model for performance evaluation of higher educational institutions in Brazil. In this study, they proposed an approach that is used to measure the efficiencies of 45 Brazilian federal universities. Yang et al. \cite{yang2018measuring} also applied a two-stage network DEA model for analyzing the reasons behind the inefficiencies of Chinese research universities. Moreno-G{\'o}mez et al. \cite{moreno2019measuring} measured efficiencies of 78 Colombian universities across the period 2015–2017 by using a two-stage network DEA model. Monfared and Safi \cite{monfared2013network}, proposed a set of performance indicators to enable efficient analysis of academic activities and apply network DEA structure to analyze the performance of the academic colleges at Alzahra University in Iran.

 There are many kinds of research available in the field of network DEA, but the application of network DEA in the field of the education sector is still very limited. In the current study, by considering quantitative attributes to measure the performance of Indian institutes of management (IIMs) by using network DEA, we develop a procedure that captures both quality and quantity. In the current study, we developed a research production model to elaborate the research process in IIMs. With the help of a two-stage relational network DEA model performance analysis is done. The results are compared with the conventional DEA model. In this study, we aim to know the reasons for the inefficiencies of DMUs with the help of the network DEA model.

The rest of the paper is organized into 4 sections. Introduction is  given in Section \ref{itm:intro}. A two-stage relational network DEA model and supported statements are given in Section \ref{itm:Twostagemodel}. The conceptual framework of the proposed two-stage research production model and education sector application with results is discussed in Section \ref{itm:appli}. Finally, some concluding remarks and future scope are given in Section \ref{itm:conclu}.

\section{ Two-stage relational network DEA model}\label{itm:Twostagemodel}

Suppose there are $n$ homogeneous DMUs ($DMU_{j};j=1, 2, 3,..., n$), which consumes $m$ inputs, $ x_{ij},\,i=1, 2, 3,..., m$ to produce $s$ outputs  $ y_{rj}, \,r=1, 2, 3,...,s$. Then the efficiency score of the $k^{th}$ DMU can be defined with the help of the following CCR DEA model\cite{cooper2007data}:

\begin{description}
	\item[{Model 1} \label{itm:Model1} ]:  For $k=1,2,3,...,n,$
\begin{flalign} \label{Model1}
	& \max E_k= \frac{\sum_{r=1}^{s}y_{rk}v_{rk}}{\sum_{i=1}^{m}x_{ik}u_{ik}} &\\
	&\mbox{subject to}~~
	\frac{\sum_{r=1}^{s}y_{rj}v_{rk}}{\sum_{i=1}^{m}x_{ij}u_{ik}} \leq 1; &\\
&u_{ik}\geq \varepsilon\,\,\forall i,  v_{rk}\geq \varepsilon\,\,\forall r,\,\,\,\varepsilon > 0;\nonumber&
\end{flalign}
where $\varepsilon$ is the non-Archimedean infinitesimal; $y_{rj}$ is the amount of the $r^{th}$ output produced by the $j^{th}$ DMU; $x_{ij}$ is the amount of the $i^{th}$ input used by the $j^{th}$ DMU; $u_{ik}$ and $v_{rk}$ are the weights corresponding to the $i^{th}$ input and $r^{th}$ output respectively.
\end{description}

\begin{description}
	\item[\textbf{Definition 1:}\cite{charnes1978measuring}\label{itm: Def 1}] 
	$E_{k}$ is called the CCR efficiency of $DMU_k$. $E_{k}$  lies in the interval (0,1]. Let $E_{k}^{*}$ be the optimal value of $E_{k}$. Then $DMU_{k}$ is said to be {\bf CCR-efficient} if $E_{k}^{*}=1$, an, $DMU_{k}$ is said to be {\bf CCR-inefficient} if $E_{k}^{*}<1$. \\
		\end{description}
		

Now, we consider a production process composed of two sub-processess as depicted in Fig \ref{fig: Fig. 1.}. The whole process consists of two stages. The 1st stage uses m inputs $x_{ik}, ~ i=1,2,...,m$  to produce p intermediate products $z_{dk},~ d=1,2,3,...,p$. The 2nd stage uses p intermediate products $z_{dk},~ d=1,2,3,...,p$  to produce s outputs $y_{rk},~ r=1,2,3,...,s$. The whole process uses m inputs $x_{ik}, ~ i=1,2,...,m$  to produce s outputs $y_{rk},~ r=1,2,3,...,s$. Moreover, the outputs of 1st stage (i.e. $z_{dk}$) are the inputs of the 2nd stage. The efficiency of 1st and 2nd stage can be calculated by the conventional CCR DEA model explained in \nameref{itm:Model1}. Suppose the efficiency of $k^{th}$ DMU in 1st stage is $E_{k}^{1}$ and in the 2nd stage is $E_{k}^{2}$. Then the 1st and 2nd stage efficiencies of $DMU_k$ can be calculated with the help of \nameref{itm:Model2} and \nameref{itm:Model 3} as follows \cite{kao2008efficiency}:

\begin{figure}[ht]
	\centering
	\includegraphics[width=0.8\linewidth]{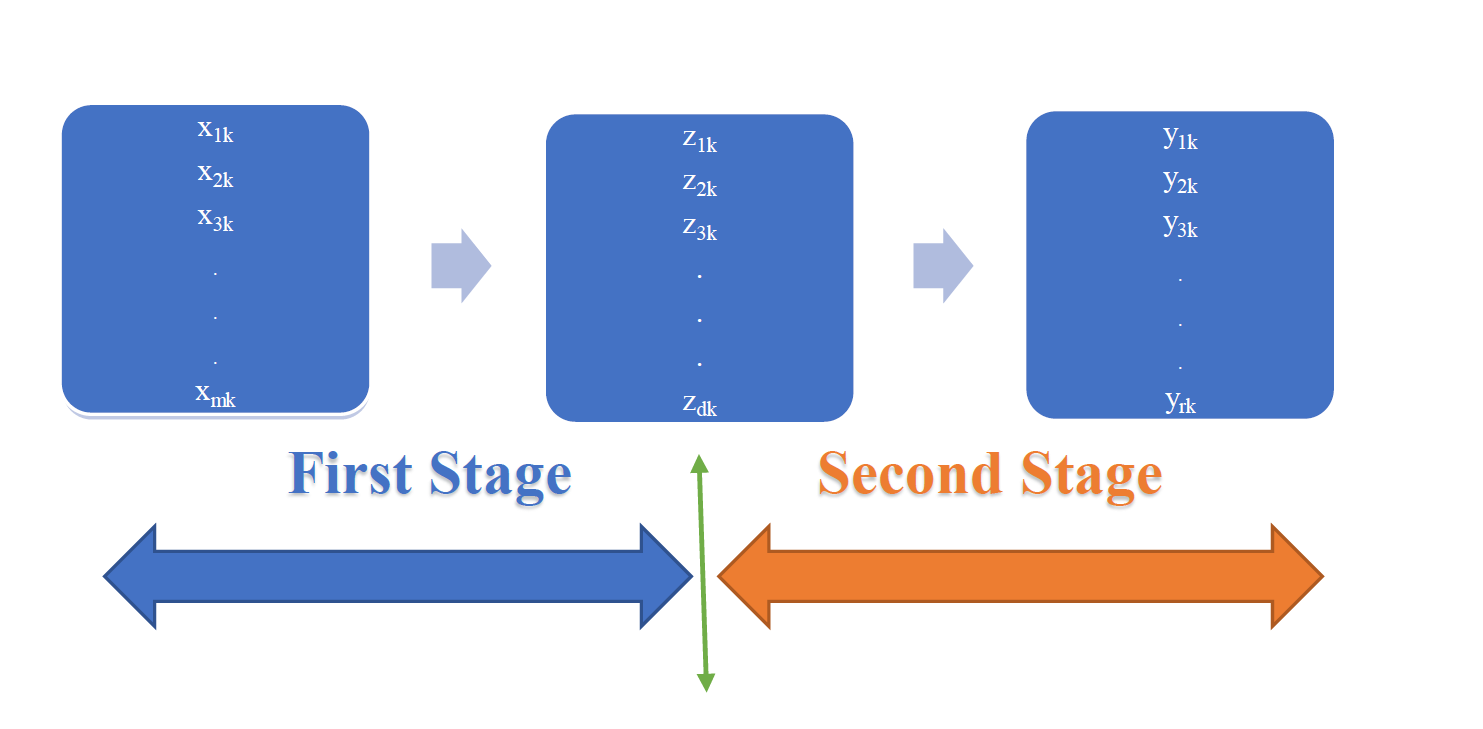} 
	\caption{A two-stage series production process system with input x, output y, and intermediate products z.}
	\label{fig: Fig. 1.}
\end{figure}

\begin{description}
	\item[{Model 2}\label{itm:Model2} ]:  For $k=1,2,3,...,n,$
	\begin{flalign} 
	& \max E_{k}^1= \frac{\sum_{d=1}^{p}z_{dk}w_{dk}}{\sum_{i=1}^{m}x_{ik}u_{ik}} &\\
	&\mbox{subject to}~~
	\frac{\sum_{d=1}^{p}z_{dj}w_{dk}}{\sum_{i=1}^{m}x_{ij}u_{ik}} \leq 1; &\\
	&u_{ik}\geq \varepsilon\,\,\forall i,  w_{dk}\geq \varepsilon\,\,\forall p,\,\,\,\varepsilon > 0.\nonumber&
	\end{flalign}
\end{description}

\begin{description}
	\item[{Model 3} \label{itm:Model 3} ]:  For $k=1,2,3,...,n,$
	\begin{flalign} 
	& \max E_{k}^2= \frac{\sum_{r=1}^{s}y_{rk}v_{rk}}{\sum_{d=1}^{p}z_{dk}w'_{dk}} &\\
	&\mbox{subject to}~~
	\frac{\sum_{r=1}^{s}y_{rj}v_{rk}}{\sum_{d=1}^{p}z_{dj}w'_{dk}} \leq 1; &\\
	&v_{rk}\geq \varepsilon\,\,\forall r,  w'_{dk}\geq \varepsilon\,\,\forall p,\,\,\,\varepsilon > 0.\nonumber&
	\end{flalign}
\end{description}

Note that the weights $w_{dk}$ and $w'_{dk}$ used in \nameref{itm:Model2} and \nameref{itm:Model 3} are same \cite{lee2016network}. The reason is, input of sub-process 2 to be the expected output of sub-process 1 while calculating the expected output of the whole process. When this concept is generalized to the case of multiple intermediate products, it requires the aggrigated value of intermediate products, which is $\sum_{d=1}^{p}w_{dk}z_{dk}$ to be the same \cite{kao2008efficiency}. Hence \nameref{itm:Model 3} can be re-written as follows: 

\begin{description}
	\item[{Model 4}\label{itm:Model4} ]:  For $k=1,2,3,...,n,$
	\begin{flalign} 
	& \max E_{k}^2= \frac{\sum_{r=1}^{s}y_{rk}v_{rk}}{\sum_{d=1}^{p}z_{dk}w_{dk}} &\\
	&\mbox{subject to}~~
	\frac{\sum_{r=1}^{s}y_{rj}v_{rk}}{\sum_{d=1}^{p}z_{dj}w_{dk}} \leq 1; &\\
	&v_{rk}\geq \varepsilon\,\,\forall r,  w_{dk}\geq \varepsilon\,\,\forall p,\,\,\,\varepsilon > 0.\nonumber&
	\end{flalign}
\end{description}

With the help of \nameref{itm:Model1}, efficiency of the whole process can be calculated independently as the 1st stage efficiency model (\nameref{itm:Model2}) and 2nd stage efficiency model (\nameref{itm:Model4}) are essentially same. (\nameref{itm:Model2}) and (\nameref{itm:Model4}) give the 1st and second stage efficiencies of independent network DEA model. To link these two sub-processes with the whole process, a model must be constructed that describe the series relationship of the whole process with the two sub-processes.\\

Suppose, for $DMU_k$; $u_{ik}^{*}$, $v_{rk}^{*}$, $w_{dk}^{*}$ are the optimal weights (multipliers) that are obtained while calculating  overall efficiency $E_k$, sub-process efficiencies $E_{k}^{1}$, and  $E_{k}^{2}$. Then we have,\\

	\begin{flalign} 	
	&	E_k=\frac{\sum_{r=1}^{s}y_{rk}v_{rk}^{*}}{\sum_{i=1}^{m}x_{ik}u_{ik}^{*}} \label{Eq: Eq:9}&\\
		&  E_{k}^1=\frac{\sum_{d=1}^{p}z_{dk}w_{dk}^*}{\sum_{i=1}^{m}x_{ik}u_{ik}^*} \label{Eq: Eq:10} &\\	
	&  E_{k}^2=\frac{\sum_{r=1}^{s}y_{rk}v_{rk}^*}{\sum_{d=1}^{p}z_{dk}w_{dk}^*}  \label{Eq: Eq:11}&
	\end{flalign}

From Eqs. \ref{Eq: Eq:9},\ref{Eq: Eq:10} \& \ref{Eq: Eq:11}, we have a relation among $E_k$, $E_{k}^{1}$, and  $E_{k}^{2}$ as;
\begin{equation}
\textbf{$E_k= E_{k}^{1} \times E_{k}^{2}$}. \label{Eq: Eq. (12)}\\
\end{equation}

Based on this concept, let us incorporate the ratio constraints of the two sub-processes into conventional DEA model i.e., \nameref{itm:Model1}, to get a two-stage relational network DEA model.  This will give the overall efficiency $E_{k}$, to a two-stage relational network DEA model, taking into account the series relationship of two sub-processes as follows:

\begin{description}
	\item[{Model 5} \label{itm:Model5} ]:  For $k=1,2,3,...,n,$
	\begin{flalign} 
	& \max E_k= \frac{\sum_{r=1}^{s}y_{rk}v_{rk}}{\sum_{i=1}^{m}x_{ik}u_{ik}} &\\
	&\mbox{subject to}~~
	\frac{\sum_{r=1}^{s}y_{rj}v_{rk}}{\sum_{i=1}^{m}x_{ij}u_{ik}} \leq 1; &\\
	&\frac{\sum_{d=1}^{p}z_{dj}w_{dk}}{\sum_{i=1}^{m}x_{ij}u_{ik}} \leq 1;&\\
	&\frac{\sum_{r=1}^{s}y_{rj}v_{rk}}{\sum_{d=1}^{p}z_{dj}w_{dk}} \leq 1; &\\
	&u_{ik}\geq \varepsilon\,\,\forall i,  v_{rk}\geq \varepsilon\,\,\forall r,\,\,\,  w_{dk}\geq \varepsilon\,\,\forall p,\varepsilon > 0.\nonumber&
	\end{flalign}
\end{description}

  Now, \nameref{itm:Model5} can be written in linear form as follows:

\begin{description}
	\item[{Model 6} \label{itm:Model6} ]:  For $k=1,2,3,...,n,$
	\begin{flalign} 
	& \max E_k= {\sum_{r=1}^{s}y_{rk}v_{rk}} &\\
	&\mbox{subject to}~~
	{\sum_{i=1}^{m}x_{ik}u_{ik}}=1;&\\
&	{\sum_{r=1}^{s}y_{rj}v_{rk}}-{\sum_{i=1}^{m}x_{ij}u_{ik}} \leq 0; &\\
	&{\sum_{d=1}^{p}z_{dj}w_{dk}}-{\sum_{i=1}^{m}x_{ij}u_{ik}} \leq 0;&\\
	&{\sum_{r=1}^{s}y_{rj}v_{rk}}-{\sum_{d=1}^{p}z_{dj}w_{dk}} \leq 0; &\\
	&u_{ik}\geq \varepsilon\,\,\forall i,  v_{rk}\geq \varepsilon\,\,\forall r,\,\,\,  w_{dk}\geq \varepsilon\,\,\forall p,\varepsilon > 0.\nonumber&
	\end{flalign}
\end{description}

It is quite possible that the optimal solution (optimal multipliers) obtained after solving \nameref{itm:Model6} are not unique. Therefore, the relation \textbf{$E_k= E_{k}^{1} \times E_{k}^{2}$} would not be unique. This makes comparison of $E_{k}^{1}$ and $E_{k}^{2}$  between all DMUs not the same basis.  One solution to this problem is to find a set of multipliers that produce the largest $E_k^{1}$ while maintaning the $E_k$ obtained from \nameref{itm:Model6}.

\begin{description}
	\item[{Model 7} \label{itm:Model7} ]:  For $k=1,2,3,...,n,$
	\begin{flalign} 
	& \max E_{k}^{1}= {\sum_{d=1}^{p}z_{dk}w_{dk}} &\\
	&\mbox{subject to}~~
	{\sum_{i=1}^{m}x_{ik}u_{ik}}=1;&\\
	&	{\sum_{r=1}^{s}y_{rk}v_{rk}}-E_{k}{\sum_{i=1}^{m}x_{ik}u_{ik}} = 0; &\\
	&	{\sum_{r=1}^{s}y_{rj}v_{rk}}-{\sum_{i=1}^{m}x_{ij}u_{ik}} \leq 0; &\\
	&{\sum_{d=1}^{p}z_{dj}w_{dk}}-{\sum_{i=1}^{m}x_{ij}u_{ik}} \leq 0;&\\
	&{\sum_{r=1}^{s}y_{rj}v_{rk}}-{\sum_{d=1}^{p}z_{dj}w_{dk}} \leq 0; &\\
	&u_{ik}\geq \varepsilon\,\,\forall i,  v_{rk}\geq \varepsilon\,\,\forall r,\,\,\,  w_{dk}\geq \varepsilon\,\,\forall p,\varepsilon > 0.\nonumber&
	\end{flalign}
\end{description}

After calculating $E_{k}^{1}$ with the help of \nameref{itm:Model7}, second stage efficiency efficiency $E_{k}^{2}$ can be calculated easily by using the relation ( Eq. \ref{Eq: Eq. (12)}) as {$ E_{k}^{2}= E_{k}/ E_{k}^{1}$}. Silmilarly, if decision-maker is more concerned for second stage efficiency then, one can formulated model for calculating $E_{k}^{2}$ by the following model as follows:

\begin{description}
	\item[{Model 8} \label{itm:Model8} ]:  For $k=1,2,3,...,n,$
	\begin{flalign} 
	& \max E_{k}^{2}= {\sum_{r=1}^{s}y_{rk}v_{rk}} &\\
	&\mbox{subject to}~~
	{\sum_{d=1}^{p}z_{dk}w_{dk}}=1;&\\
	&	{\sum_{r=1}^{s}y_{rk}v_{rk}}-E_{k}{\sum_{i=1}^{m}x_{ik}u_{ik}} = 0; &\\
	&	{\sum_{r=1}^{s}y_{rj}v_{rk}}-{\sum_{i=1}^{m}x_{ij}u_{ik}} \leq 0; &\\
	&{\sum_{d=1}^{p}z_{dj}w_{dk}}-{\sum_{i=1}^{m}x_{ij}u_{ik}} \leq 0;&\\
	&{\sum_{r=1}^{s}y_{rj}v_{rk}}-{\sum_{d=1}^{p}z_{dj}w_{dk}} \leq 0; &\\
	&u_{ik}\geq \varepsilon\,\,\forall i,  v_{rk}\geq \varepsilon\,\,\forall r,\,\,\,  w_{dk}\geq \varepsilon\,\,\forall p,\varepsilon > 0.\nonumber&
	\end{flalign}
\end{description}

The efficiency of 1st stage $E_{k}^{1}$ can be calculated easily by using the relation ( Eq. \ref{Eq: Eq. (12)}) as {$ E_{k}^{1}= E_{k}/ E_{k}^{2}$}.

%
%
%

	\begin{theorem}\label{theorem : theorem 3.1}
		For $DMU_k$, if $E_{k}=1~ \implies~E_{k}^{1}=E_{k}^{2}=1$. Where $E_{k}^{1}$ and $E_{k}^{2}$ are the optimal values of 1st and 2nd stage efficiencies obtained from \nameref{itm:Model7} and \nameref{itm:Model8}.
	\end{theorem}
	
	\begin{description}
		\item[Proof:] Suppose $(u_{ik}^{*}, v_{rk}^{*},w_{dk}^{*}~ \forall~ i,r,d)$  be an optimal solution of \nameref{itm:Model6} then we have\\
		
		$E_k= {\sum_{r=1}^{s}y_{rk}v_{rk}^*}$,
		
		\begin{description}
			\item[] and
			\begin{flalign} 
			&{\sum_{i=1}^{m}x_{ik}u_{ik}^*}=1;&\\
			&	{\sum_{r=1}^{s}y_{rj}v_{rk}^*}-{\sum_{i=1}^{m}x_{ij}u_{ik}^*} \leq 0; &\\
			&{\sum_{d=1}^{p}z_{dj}w_{dk}^*}-{\sum_{i=1}^{m}x_{ij}u_{ik}^*} \leq 0;&\\
			&{\sum_{r=1}^{s}y_{rj}v_{rk}^*}-{\sum_{d=1}^{p}z_{dj}w_{dk}^*} \leq 0; &\\
			&u_{ik}^*\geq \varepsilon\,\,\forall i,  v_{rk}^*\geq \varepsilon\,\,\forall r,\,\,\,  w_{dk}^*\geq \varepsilon\,\,\forall p,\varepsilon > 0.\nonumber&
			\end{flalign}
		\end{description}
		With $E_{k}=1$, we have, for specific $DMU_k:$\\
		\begin{equation}
		1={\sum_{r=1}^{s}y_{rk}v_{rk}^*}\leq {\sum_{d=1}^{p}z_{dk}w_{dk}^*} \leq {\sum_{i=1}^{m}x_{ik}u_{ik}^*} =1 \label{Eq: Eq (38)}
		\end{equation}
		
		Therefore, from Eq. (\ref{Eq: Eq (38)})	it is possible that ${\sum_{r=1}^{s}y_{rk}v_{rk}^*}= {\sum_{d=1}^{p}z_{dk}w_{dk}^*}=1$. In other words, the optimal solution $(u_{ik}^{*}, v_{rk}^{*},w_{dk}^{*}~ \forall~ i,r,d)$  can make the objective function values of \nameref{itm:Model7} and \nameref{itm:Model8} are equal to 1. It is easy to varify that $(u_{ik}^{*}, v_{rk}^{*},w_{dk}^{*}~ \forall~ i,r,d)$  can satisfy all the constraints of \nameref{itm:Model7} and \nameref{itm:Model8} simultaneously.\\
		Thus, we have $E_{k}^{1}=E_{k}^{2}=1$ for $DMU_k.$
		
		$~~~~~~~~~~~~~~~~~~~~~~~~~~~~~~~~~~~~~~~~~~~~~~~~~~~~~~~~~~~~~~~~~~~~~~~~~~~~~~~~~~~~~~~~~~~~~~~~~~~~~~~~~~~~~~~~~~~~~~~~~~~~~~~~~~~~~~~\square$		
	\end{description}
	
	\section{Education sector application}\label{itm:appli}
	For the application of the network DEA, we consider the problem of checking the performance efficiency of the Indian Institutes of Management (IIMs) in India. IIMs are administrative and research educational institutions in India. Primarily they offer undergraduate, postgraduate, doctoral, and other additional courses. IIMs have been declared the most important institutions in the country by The Ministry of Education Govt. of India.
	
	\subsection{Conceptual framework} 
\begin{figure}[ht]
		\centering
		\includegraphics[width=0.8\linewidth]{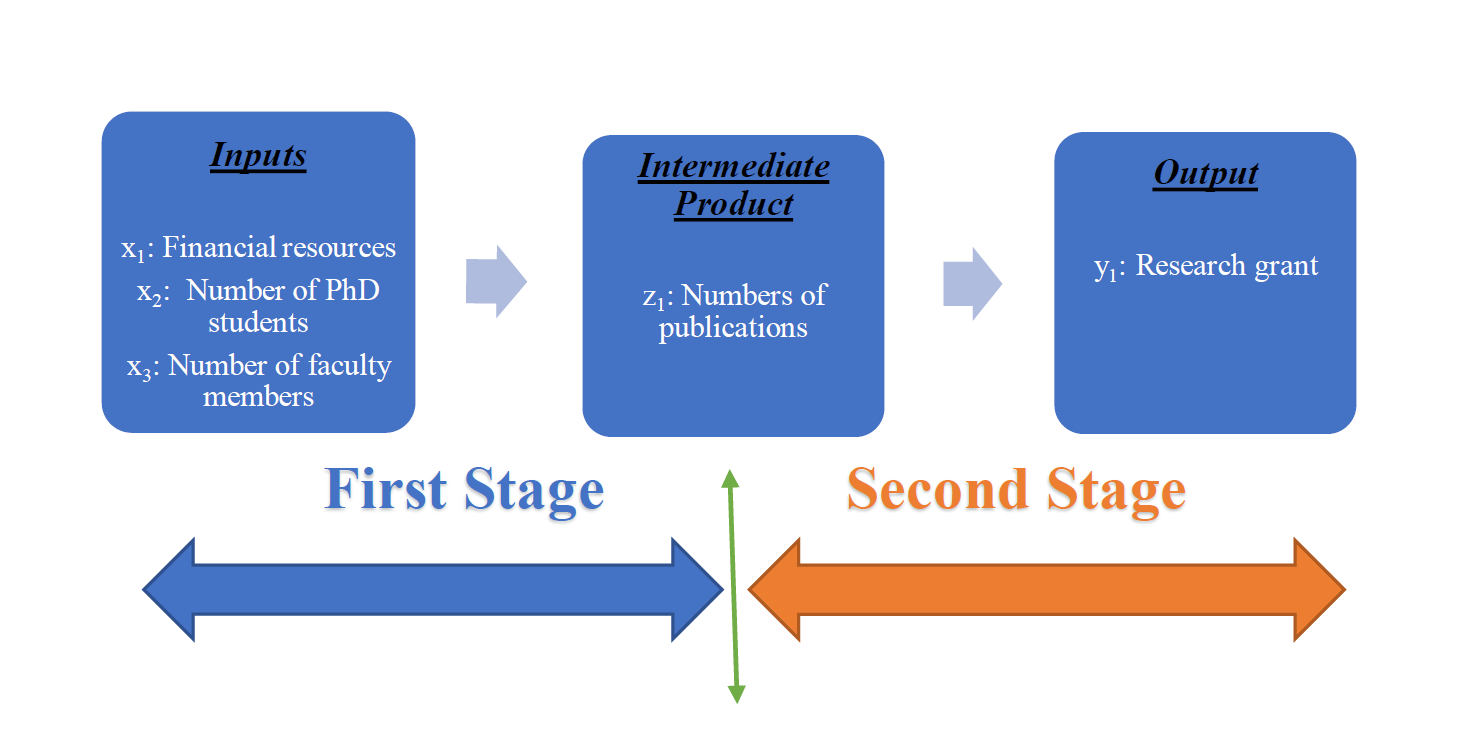} 
		\caption{A two-stage series production process system with inputs $x_1, x_2, x_3$, output y, and intermediate products z.}
		\label{fig: Fig. 2.}
	\end{figure}

	Education and research play a vital role in the development of any nation. In the current study, we will try to measure the performance efficiencies of IIMs with the help of the network DEA model, particularly the two-stage network DEA. The main reason for using the two-stage network DEA model over the conventional CCR DEA model is to understand the causes behind the inefficiencies of DMUs. The network DEA model used for the current study is depicted in Fig \ref{fig: Fig. 2.}. Here, in this study, we constructed a two-stage research production model of IIMs. The research production model in IIMs is complex. Fig \ref{fig: Fig. 2.} depicts our network DEA model for IIMs. This model consists of two stages. In the 1st stage, the model consumes $x_{ij}$ inputs to produce $z_{dj}$ intermediate products. Outputs of 1st stage start working as inputs in 2nd stage.
	
 In our research production model, we use three inputs ($x_1$: financial resources, $x_2$: number of Ph.D. students, $x_1$: number of faculty members) in 1st stage. In research, these three parameters play a very crucial role. For any Ph.D. student, it is always challenging to research without a guide and financial support. This is the reason why we considered these three parameters as inputs in 1st stage. The output of 1st stage or intermediate product is considered as the 'number of publications' ($z_1$). This is the most obvious output for any research production model. In India, Ph.D. students have to publish a few research papers to submit their thesis. These publications validate the choice of inputs as, Ph.D. students, faculties, and financial resources. The output of the 1st stage $z_1$ is used as an intermediate input in the 2nd stage to produce the final output $y_1$ (research grant). Since in the present study, we are focusing only on research activities so we are not very much bothered about other parameters. For any IIM, it is very important to win ample research grants so that they can continue their research process. IIMs get their research grant from two types of research projects; first as the sponsored research projects, and second as consultancy research projects. In our study, we consider the research grant index ($y_1$) as the sum of research grants obtained from sponsored research projects and consultancy research projects. It is a widespread perception that can be treated as the fact that if any IIMs publishes more research papers then that IIM can win more research grants. Hence the choice of intermediate input $z_1$ is the best choice for the final output $y_1$ for our two-stage research production model.\\
 
 \subsection{Results and discussion}
In the current study, IIMs are considered as our DMUs. The data is collected from the NIRF website launched by The Ministry of Education on 29th September 2015. Data is collected for the year 2020-21. Input, intermediate input, and output data for 13 IIMs for the application two-stage network DEA model is given in Table \ref{tab:Table-1}.\\

\begin{table}[ht]
	\centering
	
	\caption{Inputs (x), intermediate product (z), and output(y) for 13 IIMs}
	\label{tab:Table-1}
	\resizebox{\textwidth}{!}{
		\begin{tabular}{|llllllll|}
			\hline
			DMU & IIM Name      & State  \qquad\qquad  & \multicolumn{3}{l}{\qquad \qquad Inputs} & \multicolumn{1}{l}{Intermediate product}& \multicolumn{1}{l|}{Output} \\ \hline
			&               &           & $\qquad{x_1}$   & ${x_2}$  & ${x_3}$  &$ {\qquad \qquad z_1}$ &$\qquad {y_1}$  \\ \hline
			$D_1$  & IIM Bangalore & Karnataka&2783060307	&100	  &116  &\qquad \qquad 81	& 1302218468      \\
			
			$D_2$  & IIM Ahmedabad & Gujrat &2282349127	&128&110&\qquad \qquad 43&   34953142     \\
			
			$D_3$  & IIM Calcutta & West Bengal &1372017597&
			97	&91&\qquad \qquad 94&261107000      \\
			
			$D_4$  & IIM Lucknow & Uttar Pradesh &1405523203
				&178&104&\qquad \qquad 89&10643365
				      \\
			$D_5$  & IIM Indore & Madhya Pradesh &1323781000
				&120&139&\qquad \qquad 124&165125234
				     \\
			
			$D_6$  & IIM Kozhikode & Kerala &1009040766
				&91&109&\qquad \qquad47 &  112206310
				   \\
			
			$D_7$  & IIM Udaipur & Rajasthan &675751470
				&13&42&\qquad \qquad 27 &   1091766
				     \\
			
			$D_8$  & IIM Tiruchirapalli & Tamilnadu &366034016
				&22&55&\qquad \qquad25 & 6119078
				      \\
			
			$D_9$  & IIM Raipur & Chhatisgarh  &724040010
				&21&36&\qquad \qquad 40 &   9499507
				        \\
			
			$D_{10}$  & IIM Rohtak & Haryana &343393901
				&53&92&\qquad \qquad 47&	26337325
				       \\
			
			$D_{11}$  & IIM Shillong & Meghalaya &373246687
				&23&32&\qquad \qquad 44&    5971808
				   \\
			
			$D_{12}$  & IIM Kashipur & Uttarakhand &337690753
				&54&45&\qquad \qquad31 &    8121608
				   \\
			
			$D_{13}$  & IIM Ranchi & Jharkhand &287052918
				&66&42&\qquad \qquad 65&   696200
				   \\	
			\hline
		\end{tabular}
	}
\end{table}

       \begin{table}[ht]
       	\centering
       	
       	\caption{Efficiency scores and ranks (in parantheses) of two-stage network DEA model}
       	\label{tab:Table-2}
       	\resizebox{\textwidth}{!}{
       		\begin{tabular}{|lllll|}
       			\hline
       			DMU & IIM Name        & \multicolumn{3}{l|}{ Two-stage relational network DEA model}  \\ \hline
       			& & $~~~E_k$& $~~~E_k^1$& $~~~E_k^2$ \\ \hline
       			$D_1$  & IIM Bangalore  &0.4973(1)&0.4973(9)&~~~ 1(1)      \\
       			
       			$D_2$  & IIM Ahmedabad &0.0135(7)&0.2668(11)&0.0506(5)     \\
       			
       			$D_3$  & IIM Calcutta &0.1235(2)&0.7147(4)&0.1728(2)     \\
       			
       			$D_4$  & IIM Lucknow  &0.0041(11)&0.5529(7)&0.0074(11)
       			\\
       			$D_5$  & IIM Indore   &0.0568(3)&0.6857(5)&0.0823(4)
       			\\
       			
       			$D_6$  & IIM Kozhikode &0.0507(4)&0.3417(10)&0.1485(3)
       			\\
       			
       			$D_7$  & IIM Udaipur &0.0025(12)&1(1)&0.0025(12)
       			\\
       			
       			$D_8$  & IIM Tiruchirapalli &0.0090(8)&0.5934(6)&0.0152(8)
       			\\
       			
       			$D_9$  & IIM Raipur  &0.0145(6)&0.9809(2)&0.0148(9)
       			\\
       			
       			$D_{10}$  & IIM Rohtak &0.0268(5)&0.7703(3)&0.0348(6)
       			\\
       			
       			$D_{11}$  & IIM Shillong &0.0084(9)&1(1)&0.0084(10)
       			\\
       			
       			$D_{12}$  & IIM Kashipur  &0.0082(10)&0.5064(8)&0.0163(7)
       			\\
       			
       			$D_{13}$  & IIM Ranchi &0.0007(13)&1(1)&0.0007(13)
       			\\	
       			\hline
       		\end{tabular}
       	}
       \end{table}

 In our research production model, the efficiency of 1st stage measures the performance of IIMs based on publications, while the 2nd- stage measures the performance on the basis of the research grant won by any IIM. By applying \nameref{itm:Model6}, the overall efficiency of two-stage relational network DEA model for 13 IIMs is calculated. The results are shown in Table \ref{tab:Table-2} under the subcolumn $E_k$ of the column two-stage relational network DEA model. Since, to win an attractive research grant is of more concern for any IIM so that the research process can be performed without any obstacle, we choose to measure $E_k^{2}$ first with the help of $E_k$ and after this calculate $E_k^{1}$ via $E_k^{1}={E_k}/{E_k^{2}}$. The choice of selecting performance measures as a research grant or a number of publications is completely dependent upon the decision-maker. Here, the research grants are our primary concern, while publications can be the topic of paramount concern for some other decision-maker. This choice is not going to affect the 1st and 2nd stage efficiency scores of a two-stage relational network DEA model.\\

Now we will analyze the results obtained from applying \nameref{itm:Model6}, \nameref{itm:Model7}, and \nameref{itm:Model8} on the data given in Table \ref{tab:Table-1}. The results are given under the column named as \textbf{'two-stage relational network DEA model'}. Note that in Table \ref{tab:Table-2}, none of the 13 IIMs, performed efficiently in both stages. IIM Bangalore achieved the highest overall efficiency score of 0.4973 while, IIM Ranchi achieved the lowest overall efficiency score of 0.0007. For the 1st stage, there are 3 IIMs named, IIM Udaipur, IIM Shillong, and IIM Ranchi, which perform efficiently; as they have an efficiency score $E_k^{1}=1.$ In the 2nd stage, IIM Bangalore is the only IIM that is efficient. As we have already derived a relation \textbf{$E_k= E_{k}^{1} \times E_{k}^{2}$}  in Eq. (\ref{Eq: Eq. (12)}), indicates that the overall efficiency is the product of 1st and 2nd stage efficiencies. So there is no $E_k$ greater than its corresponding $E_{k}^{1} ~or~ E_{k}^{2}$. It is to be noted that IIM Udaipur ($D_7$), IIM Shillong ($D_{11}$), and IIM Ranchi ($D_{13}$) perform efficiently in 1st stage but not efficiently in the 2nd stage. It means that these IIMs are good at publishing research but still not able to win an attractive research grant. Whereas IIM Bangalore ($D_1$) does not perform well in the 1st stage but in the 2nd stage it is performing efficiently. This indicates that it is a bit easy for IIM Bangalore to win ample research grants. It is observable from Table \ref{tab:Table-2} that except for IIM Bangalore, the efficiency of the 1st stage is higher than the efficiency of the 2nd stage for all IIMs. This indicates that the low $E_k$ of the whole process is mainly due to the low $E_{k}^{2}$ of the 2nd stage (i.e., research grant). Since $E_k$ is always less than or equal to $E_{k}^{1}$, and $E_{k}^{2}$ so it will be more informative to look at the ranks, which are the numbers in brackets in Table \ref{tab:Table-2}, over the efficiency scores. If any IIM has similar ranks in $E_{k},~E_{k}^{1}, ~and~ E_{k}^{2}$ then it can be said that the particular IIM performs evenly in the whole process and in both the sub-processes. Hence, with the help of a two-stage relational network DEA, we are in the position to understand the reason behind the inefficiency of any IIM.\\

Kao and Hung \cite{kao2008efficiency} explained logically that $E_k$ for any IIM should lie either between $E_{k}^{1}, ~and~ E_{k}^{2}$ or in the neighborhood of $E_{k}^{1}, ~and~ E_{k}^{2}$. The main reason for this is that the performance of the whole process is the aggregation  of two sub-processes. Here in Table \ref{tab:Table-2}, it can be seen that out of 13 IIMs the ranks of 10 IIMs lies between $E_{k}^{1}, ~and~ E_{k}^{2}$, and the ranks of remaining 3 IIMs lies in the neighborhood of $E_{k}^{1}, ~or~ E_{k}^{2}$.\\

Now, with the help of two-stage relational network DEA model we are  some what able to understand the actual reasons of inefficiencies in DMUs.  As we know that the DEA works like a black box that does not provide any adequate detail to identify the specific reason for inefficiency in decision-making units (DMUs). After knowing the reasons of inefficiecies, now we are interested in the comparison of the efficiency of the whole research production process with the two-stage relational network DEA model ($E_k$) and conventional CCR DEA model ($E_k^{CCR}$). Actually, we wanted to understand what will happen if the overall efficiency of DMUs is calculated by taking sub-stages of the model in mind and without bothering about sub-stages. The comparison is done and the results are shown in Table \ref{tab:Table-3}. The results depicted in Table \ref{tab:Table-3} shows that the conventional CCR efficiency ($E_k^{CCR}$) is higher than two-stage relational network DEA efficiency ($E_k$). But it will be more appropriate to compare the ranks of IIMs rather than their efficiency scores because the relational model has two more sets of constraints than the conventional CCR DEA model. These additional sets of constraints represent the 1st and 2nd stage efficiency scores. The 1st and 2nd stage efficiency score can be affecting factors in the overall efficiency scores of IIMs in the two-stage relational network DEA model as explained earlier. Kao and Hung \cite{kao2008efficiency}, mathematically, explained in detail why this deviation occurs in the $E_k$ and $E_k^{CCR}$ by constructing the dual model of the relational two-stage network DEA model. A Spearman's rank \cite{zar2005spearman} correlation coefficient is calculated for $E_k$ and $E_k^{CCR}$. The correlation coefficient $(\rho=0.91758)$ indicates that efficiencies obtained from both two-stage relational and conventional CCR DEA models are highly correlated. This means that the efficiency scores from both models produce almost the same rankings. This can be easily observed in Table \ref{tab:Table-3}.

\begin{table}[ht]
	\centering
	
	\caption{Overall efficiency scores and ranks (in parantheses) of two-stage network DEA model and conventional CCR DEA model}
	\label{tab:Table-3}
	\resizebox{\textwidth}{!}{
		\begin{tabular}{|llll|}
			\hline
			DMU & IIM Name        & \multicolumn{1}{l}{ Two-stage relational network DEA model} & \multicolumn{1}{l|}{ Conventional CCR DEA model} \\ \hline
			& & $~~~~~~~~~~~~~~~~~E_k$& $~~~~~~~~~~~~~~~E_k^{CCR}$\\ \hline
			$D_1$  & IIM Bangalore  &~~~~~~~~~~~~~0.4973(1)&~~~~~~~~~~~~~1(1)      \\
			
			$D_2$  & IIM Ahmedabad &~~~~~~~~~~~~~0.01349(7)&~~~~~~~~~~~~~0.0327(10)     \\
			
			$D_3$  & IIM Calcutta &~~~~~~~~~~~~~0.1235(2)&~~~~~~~~~~~~~0.4067(2)     \\
			
			$D_4$  & IIM Lucknow  &~~~~~~~~~~~~~0.0041(11)&~~~~~~~~~~~~~0.0162(11)
			\\
			$D_5$  & IIM Indore   &~~~~~~~~~~~~~0.0568(3)&~~~~~~~~~~~~~0.2666(3)
			\\
			
			$D_6$  & IIM Kozhikode &~~~~~~~~~~~~~0.0507(4)&~~~~~~~~~~~~~0.2377(4)
			\\
			
			$D_7$  & IIM Udaipur &~~~~~~~~~~~~~0.0025(12)&~~~~~~~~~~~~~0.0064(12)
			\\
			
			$D_8$  & IIM Tiruchirapalli &~~~~~~~~~~~~~0.0090(8)&~~~~~~~~~~~~~0.0357(7)
			\\
			
			$D_9$  & IIM Raipur  &~~~~~~~~~~~~~0.0145(6)&~~~~~~~~~~~~~0.0347(8)
			\\
			
			$D_{10}$  & IIM Rohtak &~~~~~~~~~~~~~0.0268(5)&~~~~~~~~~~~~~0.1639(5)
			\\
			
			$D_{11}$  & IIM Shillong &~~~~~~~~~~~~~0.0084(9)&~~~~~~~~~~~~~0.0342(9)
			\\
			
			$D_{12}$  & IIM Kashipur  &~~~~~~~~~~~~~0.0082(10)&~~~~~~~~~~~~~0.0514(6)
			\\
			
			$D_{13}$  & IIM Ranchi &~~~~~~~~~~~~~0.0007(13)&~~~~~~~~~~~~~0052(13)
			\\	
			\hline
		\end{tabular}
	}
\end{table}   

\section{Conclusions}    \label{itm:conclu}                 

The main objective of the current study is to illustrate the process of research of IIMs with the help of the proposed research production model. This study applies a two-stage relational network DEA model for the performance analysis of IIMs. In our proposed research production model there are two sub-processes. As we know that higher educational institutions are complex structures that use multiple inputs to produce multiple outputs. The research process is also complicated. As, in some stage, some output of research process become the inputs for the next stages to produce final outputs. Our proposed research production model explains the use of inputs and outputs in sub-stages very interestingly. Hence, it is pretty reasonable to consider the research production model for  IIM in sub-stages, not as simply one that transforms input to produce outputs. As we know that in general, data envelopment analysis (DEA) works like a black box that does not provide any adequate detail to identify the specific reason for inefficiency in decision-making units (DMUs). 

Our study shows that the efficiency level of IIMs is higher in the conventional DEA model than the two-stage relational network DEA model. The reason for the deviation in both types of efficiency levels is explained in our study. So it will be more appropriate to compare the ranks of IIMs rather than their efficiency scores. A Spearman's rank  correlation coefficient is calculated for $E_k$ and $E_k^{CCR}$. The correlation coefficient $(\rho=0.91758)$ indicates that efficiencies obtained from both two-stage relational and conventional CCR DEA models are highly correlated. This means that the efficiency scores from both models produce almost the same rankings. But due to the use of the two-stage relational network DEA model, we are now able to know the reasons behind the inefficiencies of IIMs. We found that many IIMs performed well on the 1st stage but not on the 2nd stage. Based on this, we can conclude that for IIMs it is easier to produce research publications than to win attractive research grants. In a real-world situation, input-output data can be imprecise. The future research direction can lead us to the development and application of a two-stage network DEA model in a fuzzy environment.

\bibliographystyle{ieeetr}
\bibliography{reference}

\begin{thebibliography}{10}

\bibitem{cook2009data}
W.~D. Cook and L.~M. Seiford, ``Data envelopment analysis (dea)--thirty years
  on,'' {\em European journal of operational research}, vol.~192, no.~1,
  pp.~1--17, 2009.

\bibitem{charnes1978measuring}
A.~Charnes, W.~W. Cooper, and E.~Rhodes, ``Measuring the efficiency of decision
  making units,'' {\em European journal of operational research}, vol.~2,
  no.~6, pp.~429--444, 1978.

\bibitem{kao2014network}
C.~Kao, ``Network data envelopment analysis: A review,'' {\em European journal
  of operational research}, vol.~239, no.~1, pp.~1--16, 2014.

\bibitem{tan2021investigation}
Y.~Tan and D.~Despotis, ``Investigation of efficiency in the uk hotel industry:
  A network data envelopment analysis approach,'' {\em International Journal of
  Contemporary Hospitality Management}, 2021.

\bibitem{kao2008efficiency}
C.~Kao and S.-N. Hwang, ``Efficiency decomposition in two-stage data
  envelopment analysis: An application to non-life insurance companies in
  taiwan,'' {\em European journal of operational research}, vol.~185, no.~1,
  pp.~418--429, 2008.

\bibitem{arya2021development}
A.~Arya and S.~Singh, ``Development of two-stage parallel-series system with
  fuzzy data: A fuzzy dea approach,'' {\em Soft Computing}, vol.~25, no.~4,
  pp.~3225--3245, 2021.

\bibitem{ma2018evaluating}
J.~Ma and L.~Chen, ``Evaluating operation and coordination efficiencies of
  parallel-series two-stage system: A data envelopment analysis approach,''
  {\em Expert Systems with Applications}, vol.~91, pp.~1--11, 2018.

\bibitem{puri2017new}
J.~Puri, S.~P. Yadav, and H.~Garg, ``A new multi-component dea approach using
  common set of weights methodology and imprecise data: an application to
  public sector banks in india with undesirable and shared resources,'' {\em
  Annals of Operations Research}, vol.~259, no.~1, pp.~351--388, 2017.

\bibitem{lotfi2013multi}
F.~H. Lotfi and M.~Vaez-Ghasemi, ``Multi-component efficiency with shared
  resources in commercial banks,'' {\em International Journal of Applied},
  vol.~3, no.~4, pp.~93--104, 2013.

\bibitem{singh2021performance}
A.~P. Singh, S.~P. Yadav, and P.~Tyagi, ``Performance assessment of higher
  educational institutions in india using data envelopment analysis and
  re-evaluation of nirf rankings,'' {\em International Journal of System
  Assurance Engineering and Management}, pp.~1--12, 2021.

\bibitem{arcelus1997efficiency}
F.~ARCELUS and D.~COLEMAN, ``{An efficiency review of university
  departments},'' {\em International journal of systems science}, vol.~28,
  no.~7, pp.~721--729, 1997.

\bibitem{sinuany1994academic}
Z.~Sinuany-Stern, A.~Mehrez, and A.~Barboy, ``{Academic departments efficiency
  via DEA},'' {\em Computers \& Operations Research}, vol.~21, no.~5,
  pp.~543--556, 1994.

\bibitem{bessent1982application}
A.~Bessent, W.~Bessent, J.~Kennington, and B.~Reagan, ``{An application of
  mathematical programming to assess productivity in the Houston independent
  school district},'' {\em Management Science}, vol.~28, no.~12,
  pp.~1355--1367, 1982.

\bibitem{tomkins1988experiment}
C.~Tomkins and R.~Green, ``{An experiment in the use of data envelopment
  analysis for evaluating the efficiency of UK university departments of
  accounting},'' {\em Financial Accountability \& Management}, vol.~4, no.~2,
  pp.~147--164, 1988.

\bibitem{tyagi2009relative}
P.~Tyagi, S.~P. Yadav, and S.~Singh, ``Relative performance of academic
  departments using dea with sensitivity analysis,'' {\em Evaluation and
  Program Planning}, vol.~32, no.~2, pp.~168--177, 2009.

\bibitem{lee2016network}
B.~L. Lee and A.~C. Worthington, ``A network dea quantity and
  quality-orientated production model: An application to australian university
  research services,'' {\em Omega}, vol.~60, pp.~26--33, 2016.

\bibitem{tavares2021proposed}
R.~S. Tavares, L.~Angulo-Meza, and A.~P. Sant'Anna, ``A proposed multistage
  evaluation approach for higher education institutions based on network data
  envelopment analysis: A brazilian experience,'' {\em Evaluation and Program
  Planning}, vol.~89, p.~101984, 2021.

\bibitem{yang2018measuring}
G.-l. Yang, H.~Fukuyama, and Y.-y. Song, ``{Measuring the inefficiency of
  Chinese research universities based on a two-stage network DEA model},'' {\em
  Journal of Informetrics}, vol.~12, no.~1, pp.~10--30, 2018.

\bibitem{moreno2019measuring}
J.~Moreno-G{\'o}mez, J.~Calleja-Blanco, and G.~Moreno-G{\'o}mez, ``Measuring
  the efficiency of the colombian higher education system: a two-stage
  approach,'' {\em International Journal of Educational Management}, 2019.

\bibitem{monfared2013network}
M.~A.~S. Monfared and M.~Safi, ``Network dea: an application to analysis of
  academic performance,'' {\em Journal of Industrial Engineering
  International}, vol.~9, no.~1, pp.~1--10, 2013.

\bibitem{cooper2007data}
W.~W. Cooper, L.~M. Seiford, and K.~Tone, {\em Data envelopment analysis: a
  comprehensive text with models, applications, references and DEA-solver
  software}, vol.~2.
\newblock Springer, 2007.

\bibitem{zar2005spearman}
J.~H. Zar, ``Spearman rank correlation,'' {\em Encyclopedia of biostatistics},
  vol.~7, 2005.

\end{thebibliography}
\newpage
\end{document}